\documentstyle[11pt,leqno,rotate,epsfig]{article}
\input amssym.def
\title{A comparison of the eigenvalues of the Dirac and Laplace operator on the two-dimensional torus. \footnote {Supported by the SFB 288 of the DFG.}}
\author{ILKA AGRICOLA, BERND AMMANN and THOMAS FRIEDRICH \\
}
\date{\today}

\newcommand{\D}{\displaystyle}
\newcommand{\upsp}{\phantom{l}}
\newcommand{\downsp}{\phantom{q}}

\begin{document}

\maketitle

\mbox{} \hrulefill \mbox{}\\

\newcommand{\vol}{\mbox{vol} \, }
\newcommand{\grad}{\mbox{grad} \, }

\begin{abstract}We compare the eigenvalues of the Dirac and Laplace operator on a two-dimensional torus with respect to the trivial spin structure. In particular, we compute their variation up to order 4 upon deformation of the flat metric, study the corresponding Hamiltonian and discuss several families of examples. \end{abstract}

\vspace{0.5cm}

{\it Subj. Class.:} Differential geometry.

\mbox{} \hrulefill \mbox{}\\
\section{Introduction}

\newfont{\graf}{eufm10}
\newcommand{\alth}{\mbox{\graf h}}

We consider a two-dimensional torus $T^2$ equipped with a flat metric $g_o$ as well as a conformally equivalent metric $g$,

\[ g= h^4 g_o . \]

Denote by $\Delta_g$ the Laplace operator acting on functions and let $D_g$ be the Dirac operator. The following estimates for the first positive eigenvalues $\mu_1 (g)$ and $\lambda_1^2 (g)$ of $\Delta_g$ and $D^2_g$ are known (see [2] and [6]): 

\begin{itemize}
\item[a.)] $\D \frac{\lambda_1^2 (g_o)}{h^4_{\max}} \le \lambda_1^2 (g) \le  \frac{\lambda_1^2 (g_o)}{h^4_{\min}} \quad , \quad \frac{\mu_1 (g_o)}{h^4_{\max}} \le \mu_1 (g) \le  \frac{\mu_1 (g_o)}{h^4_{\min}} , $
\end{itemize}

where $h_{\min}$ $(h_{\max})$ denotes the minimum (maximum) of the conformal factor.\\

\begin{itemize}
\item[b.)] $\D \mu_1 (g) \le \frac{16 \pi}{\vol (T^2, g)} . $
\end{itemize}

In case the spin structure of the torus $T^2$ is nontrivial,  the Dirac operator has no kernel and,  moreover,  there exists a constant $C$ depending on the conformal  structure fixed on $T^2$ and on the spin structure such that

\[ \lambda_1^2 (g) \ge \frac{C}{\vol (T^2,g)} \]

(see [7]). However, explicit formulas for the constants are not known. In this respect, the situation on $T^2$ clearly differs from the case of the two-dimensional sphere $S^2$, where

\[ \lambda_1^2 (g) \ge \frac{4 \pi}{\vol (S^2,g)} \]

holds for any metric $g$ (see [3], [7]). \\

In this paper we compare $\mu_1 (g)$ and $\lambda_1^2 (g)$ for the trivial spin structure and a metric with $S^1$-symmetry. We will construct deformations $g_E$ of the flat metric such that $ \vol (g_E) \equiv \vol (g_o)$ and $\mu_1 (g_E) < \lambda_1^2 (g_E)$ holds for any parameter $E \not= 0$ near zero. For this purpose we calculate, in complete generality,  the formulas for the first and second variation of the spectral functions  $\mu_1$ and $\lambda_1^2$. It turns out that for any local deformation of the flat metric,  the first minimal eigenvalue of the Laplace operator is always smaller than the corresponding eigenvalue of the Dirac operator up to second order. The question whether or not there exists a Riemannian metric $g$ on the two-dimensional torus such that $\lambda_1^2 (g) < \mu_1 (g)$ remains open. Denote by $\lambda_1^2 (g;l)$ and $\mu_1 (g;k)$ the first eigenvalue of the Dirac and Laplace operator, respectively, such that its eigenspace contains the representation of weight $l$ respectively $k$. The discussion in the final part of this paper suggests the conjecture that for same index $l=k$, the eigenvalues $\lambda_1^2 (g,l)$ and $\mu_1 (g,l)$, are closely related, more precisely, that the Laplace eigenvalue is always smaller than the Dirac eigenvalue and that their difference should be measurable by some other geometric quantity.\\

The Dirac equation for eigenspinors of index $l=0$ can be integrated explicitly. In case of index $l \not= 0$, the Hamiltonian describing  the Dirac equation is a positive Sturm-Liouville operator. First, this observation yields an upper bound for $\lambda_1^2 (g,l)$. On the other hand, it proves the existence of many eigenspinors without zeros.\\

We furthermore apply the general variation formulas in order to study the eigenvalues of the Laplace and Dirac operator for the family of metrics

\[ g_E =( 1 + E \cos (2 \pi Nt) )(dt^2 +dy^2)  \]

\bigskip

in more detail. In case $N=2$, the Laplace equation is reduced to the classical Mathieu equation. A similar reduction of the Dirac equation yields a special Sturm-Liouville equation whose solutions we shall therefore call {\em Mathieu spinors}. We investigate the eigenvalues of this equation and compute (for topological index $l=1$) the first terms in the Fourier expansion of these Mathieu spinors. These computer calculations have been done by Heike Pahlisch and our grateful thanks are due to her for this. Furthermore, we thank M.~Shubin for interesting discussions on Sturm-Liouville equations.\\ 

\section{The first positive eigenvalue of the Dirac operator and Laplace operator}

Let $g$ and $g_o$ be two conformally equivalent metrics on $T^2$. Then the Laplace and Dirac operators are related by the well-known formulas

\[ \Delta_g = \frac{1}{h^4} \, \, \Delta_o  , \]
\[ D_g = \frac{1}{h^2} \, D_o + \frac{\grad (h)}{h^3}  , \]

where $\grad (h)$ denotes the gradient of the function $h$ with respect to the metric $g_o$. Let us fix the trivial spin structure on $T^2$. In this case the kernel of the operator $D_o$ coincides with the space of all parallel spinor fields, in particular, any solution of the equation $D_o (\psi_o) =0$ has constant length. The kernel of the operator $D_g$ is given by

\[ \ker (D_g) = \left\{ \frac{1}{h} \psi_o : \, \, \, D_o (\psi_o) = 0 \right\} . \]

The square $D_g^2$ of the Dirac operator preserves the decomposition $S=S^+ \oplus S^-$ of the spinor bundle $S$. Moreover, $D_g^2$ acts on the space of all sections of $S^{\pm}$ with the same eigenvalues. Therefore, the first positive eigenvalue $\lambda_1^2 (g)$ of the operator $D_g^2$ can be computed using sections in the bundle $S^+$ only. Any section $\psi \in \Gamma (S^+)$ is given by a function $f$ and a parallel spinor field $\psi_o \in \Gamma (S^+)$:

\[  \psi =(f \cdot h) \psi_o . \]

The spinor field $\psi$ is $L^2$-orthogonal to the kernel of the operator $D_g^2$ if and only if

\[ \int\limits_{T^2} ( \psi, {\textstyle{\frac{1}{h}}} \psi_o ) dT^2_g = |\psi_o|^2 \int\limits_{T^2} fh^4 dT^2_o =0 \]

holds, where $dT^2_o$ and $dT^2_g =h^4 dT^2_o$ are the volume forms of the metrics $g_o$ and $g$. The Rayleigh quotient  for the operator $D_g^2$ is given by

\[ \frac{ \D \int\limits_{T^2 \downsp}  |D_g (\psi)|^2 dT^2_g}{ \D \int\limits_{T^2}^{\upsp}  |\psi|^2 dT^2_g} \, \, = \, \, \frac{ \D \int\limits_{T^2 \downsp }  |h \cdot \grad (f) + 2 f \grad (h)|^2 dT^2_o}{ \D \int\limits_{T^2}^{\upsp}  f^2 h^6 dT^2_o} \, \, . \]

Finally, we obtain the following formulas for the first positive eigenvalue $\mu_1 (g), \lambda_1^2 (g)$ of the Laplace operator $\Delta_g$ and the square $D_g^2$ of the Dirac operator with respect to the trivial spin structure:

\[ \mu_1 (g) = \inf \left\{ \frac{\D \int\limits_{T^2 \downsp} |\grad (f) |^2 dT^2_o}{\D \int\limits_{T^2}^{\upsp} f^2 h^4 dT^2_o} : \int\limits_{T^2} fh^4 dT^2_o =0 \right\} \]

\[ \lambda_1^2  (g) = \inf \left\{ \frac{\D \int\limits_{T^2 \downsp} | h \cdot \grad (f) + 2 f \cdot \grad (h) |^2 dT^2_o}{\D \int\limits_{T^2}^{\upsp} f^2 h^6 dT^2_o} : \int\limits_{T^2} fh^4 dT^2_o =0 \right\}  . \]

A direct calculation yields the formula 

\[ \int\limits_{T^2} | h \cdot  \grad (f) + 2 f \grad (h)|^2 dT^2_o =  \int\limits_{T^2} h^2 f \Delta_o (f) dT^2_o + \int\limits_{T^2} ( 4 | \grad (h)|^2 + \textstyle{\frac{1}{2}} \Delta_o (h^2) ) f^2 dT^2_o . \]

Let us use this formula in case that $f$ is an eigenfunction of the Laplace operator, i.e., 

\[ \Delta_o (f) = \mu_1 (g) h^4 f . \]

Then it implies the following inequality between the first eigenvalues of the Laplace and Dirac operator:

\[ \lambda_1^2 (g) \le  \mu_1 (g) + \frac{\D \int\limits_{T^2 \downsp} \left( 4 |\grad (h)|^2 + \textstyle{\frac{1}{2}} \Delta_o (h^2) \right) f^2 dT^2_o}{\D \int\limits_{T^2}^{\upsp} f^2h^6 dT^2_o } \quad . \]

We are now looking  for $L^2$-estimates in case the metric $g$ admits  an $S^1$-symmetry. Indeed, let us suppose that the metric $g$ is defined on $[0,1] \times [0,1]$ by 

\[ g = h^4 (t) g_o = h^4 (t) (dt^2 + dy^2)    , \]

where the conformal factor $h^4$ depends on the variable $t$ only. Moreover, we assume that the function ${h} \, (t)$ has the symmetry

\[ h(t) = h   (1-t) .  \]

Then any function $f(t)$ with $f(t)=-f (1-t)$ satisfies the condition

\[ \int\limits_{T^2} fh^4 dT^2_o =0 \]

and, consequently, yields upper bounds for $\mu_1 (g)$ and $\lambda_1 (g)$:

\[ \mu_1 (g) \le \frac{\D \int\limits^1_{0 \downsp} |f' (t)|^2 dt}{\D \int\limits^{1^{\upsp}}_0 f^2 (t) h^4 (t) dt} := B^u_L (g;f) \]

\[ \lambda^2_1 (g) \le \frac{\D \int\limits^1_{0 \downsp} \Big(h (t)f' (t) + 2 f (t) h' (t)\Big)^2 dt}{\D \int\limits^{1^{\upsp}}_0 f^2 (t) h^6 (t) dt} := B^u_D (g;f)  \, \, . \]

\section{The first and second variation of \boldmath ${\mu_1 (g)}$ \unboldmath and \boldmath $\lambda_1^2 (g)$ \unboldmath}

We consider a Riemannian metric

\[ g = h^4 (t) g_o = h^4 (t) (dt^2 +dy^2) \]

on $T^2$ $(0 \le t \le 1, \, \, 0 \le y \le 1)$ and denote by $E(\mu_1 (g))$ and $E(\lambda_1^2 (g))$ the eigenspaces of the Laplace operator and the Dirac operator corresponding to the first positive eigenvalue. The isometry group $S^1$ acts on these eigenspaces and therefore they decompose into irreducible representations

\[ E(\mu_1 (g))= \textstyle{\sum} (k_1) \oplus \cdots \oplus \sum (k_m) \quad \mbox{and} \quad  E(\lambda_1^2 (g))= \textstyle{\sum} (l_1) \oplus \cdots \oplus \sum (l_n) , \]

where $\sum (k)$ denotes the 1-dimensional $S^1$-representation of weight $k$.\\

{\bf Proposition 1:} {\it The weights $k_{\alpha}^2$ of the first positive eigenvalue $\mu_1 (g)$ of the Laplace operator are always bounded by one:

\[ k_{\alpha}^2 \le 1 . \]

\bigskip

The weights $l_{\beta}^2$ of the first positive eigenvalue $\lambda_1^2 (g)$ of the Dirac operator are bounded by one under the condition}

\[ \max \left( \left| \frac{h' (t)}{h(t)} \right| : 0 \le t \le 1 \right) \le 3 \pi . \]

\bigskip

{\bf Proof:} Suppose that

\[ f (t,y)= A(t) e^{2 \pi k_{\alpha} iy} \]

is an eigenfunction of the Laplace operator, $\Delta_g f= \mu_1 (g)f$. Then the function $A(t)$ is a solution of the Sturm-Liouville equation

\[ - A'' (t) = \left\{ \mu_1 (g) h^4 (t) - 4 \pi^2 k_{\alpha}^2 \right\} A(t) . \]

Then consider the function

\[ F(t,y)=A(t) e^{2 \pi iy} \]

and remark that

\[ \Delta_g F= \mu_1 (g) F + 4 \pi^2 (1- k_{\alpha}^2) \frac{F}{h^4} . \]

Since $\D \int\limits_{T^2} F dT^2_g=0$,  we obtain, in case of the first positive eigenvalue, that

\[ \mu_1 \le \frac{\D \int\limits_{T^2 \downsp} \Delta_g (F)\bar{F} dT^2_g}{\D \int\limits_{T^2}^{\upsp} |F|^2 dT^2_g} = \mu_1 + 4 \pi^2 (1-k^2_{\alpha}) \, \, \frac{\D \int\limits_{T^2 \downsp} |F|^2 dT^2_o}{\D \int\limits_{T^2}^{\upsp} |F|^2 h^4 dT^2_o } \, \, . \]

The latter inequality yields $k^2_{\alpha} \le 1$ immediately. The corresponding result for the Dirac operator follows from the formula

\[ \lambda_1^2 \le \lambda_1^2 + \frac{4 \pi^2}{\D \int\limits_{T^2}^{\upsp} |F|^2 h^4 dT^2_o} \left\{ (1-l^2) \int\limits_{T^2} |F|^2 dT^2_o + \frac{l-1}{\pi} \int\limits^1_0 |F (t)|^2 \frac{h' (t)}{h(t)} dt \right\}  , \]

where we have already used the differential equation $(**)$ for $A$ that will be derived in the next paragraph. \hfill $\Box$\\

Solutions of the Laplace equation $\Delta_g f= \mu_1 (g)f$ are given by solutions of the Sturm-Liouville equation\\

\mbox{} \hfill $\displaystyle{- A'' (t) = \{ \mu_1 (g) h^4 (t) - 4 \pi^2 k^2 \} A(t)}$ \hfill $(*)$\\

with $k=0, \pm 1$. In a similar way we can reduce the Dirac equation to an ordinary differential equation. The Dirac operator $D_g$ acts on spinor fields via the formula

\[ D_g = \frac{1}{h^2(t)} \left( \begin{array}{cc} 0&i\\i&0 \end{array} \right) \partial_t + \frac{h' (t)}{h^3 (t)} \left( \begin{array}{cc} 0&i\\i&0€\end{array} \right) + \frac{1}{h^2 (t)} \left( \begin{array}{cc} 0&-1\\1&0 \end{array} \right) \partial_y  \, \, . \]

\bigskip

Suppose that a spinor field $\psi \in \Gamma (S^+)$ is a solution of the equation $D^2_g (\psi) = \lambda_1^2 (g) \psi$. Then $\psi$ is given by a solution of the Sturm-Liouville equation\\

\mbox{} \hfill $\displaystyle{- A''(t) = \left\{ \lambda_1^2 (g) h^4 (t) + \frac{h(t) h'' (t) - 2 (h' (t))^2}{h^2 (t)} - 4 \pi^2 l^2  + 4 \pi l \, \, \frac{h'(t)}{h(t)} \right\} A(t)}$ \hfill $(**)$\\

\bigskip

with $l=0, \pm 1$. In case $l=0$, this equation can be solved.\\

{\bf Proposition 2:} {\it The eigenvalues of the Sturm-Liouville equation $(**)$ for $l=0$ are given by $(n \in {\Bbb Z})$}

\[ \lambda^2 = \frac{4 \pi^2 n^2}{\left( \D \int\limits^1_0 h^2 (t) dt \right)^2}  . \]

{\bf Proof:} The Sturm-Liouville operator

\[ H = - \frac{1}{h^4 (t)} \frac{d^2}{dt^2} - \frac{h(t)h''(t) - 2 (h'(t))^2}{h^6 (t)} \]

admits a square root, namely

\[ \sqrt{H} ( -)= \frac{i}{h^3 (t)} \frac{d}{d t} (h(t) \, - ) . \]

\bigskip

Since we have $\frac{d}{dt} h^2 (t)A(t) \bar{A} (t)=0$, any solution of the equation

\[ \frac{i}{h^3 (t)} \frac{d}{d t} (h(t) A(t))= \lambda A(t) \]

satisfies the condition

\[ |A(t)| = \frac{\mbox{const}}{h(t)} . \]

Consequently, it makes sense to define a  function $f: {\Bbb R}^1 \to {\Bbb R}^1$ by the formula

\[ h(t) A(t) = e^{if(t)} , \]

for which we easily obtain the differential equation \,  $ f' (t) = \lambda h^2 (t)$. But since $A(t)$ is a periodic solution, we have the condition

\[ 2 \pi n= \int\limits^1_0 f' (t) dt =  \lambda \int\limits^1_0 h^2 (t)dt \]

for some integer $n \in {\Bbb Z}$, thus yielding the result.\\

{\bf Corollary:} {\it Let $\lambda^2 (g)$ be an eigenvalue of the Dirac operator on the two-dimensional torus $T^2$ with respect to the trivial spin structure and a Riemannian metric

\[ g= h^4 (t)(dt^2+dy^2) \]

with isometry group $S^1$. Moreover, suppose that the eigenspinor is $S^1$-invariant $(l=0)$. Then}

\[ \lambda^2 (g) \vol (T^2,g) \ge 4 \pi^2 \]

{\it holds.}\\

{\bf Proof:} Since the volume is given by $\vol (T^2,g)= \D \int\limits^1_0 h^4 (t)dt$,  the inequality follows directly from the Cauchy-Schwarz inequality $\left( \D \int\limits^1_0 h^2(t)dt \right)^2 \le \D \int\limits^1_0 h^4 (t)dt$ \,  and the previous Proposition. \hfill $\Box$\\

{\bf Remark:} This corollary \ should be compared with the following fact. Fix a re\-presentation $\Sigma (k)$ and denote by $\mu_1 (g;k)$ the first eigenvalue of the Laplace operator such that its eigenspace contains the representation $\Sigma (k)$. In case $k \not= 0$ the solution $A(t)$ of equation $(*)$ is positive (see [8], page 207) and consequently the inequality

\[ \int\limits^1_0 \Big(4 \pi k^2 - \mu_1 (g;k)h^4 (t) \Big) dt \ge 0 \]

is valid. We thus obtain the estimate

\[ 4 \pi^2 k^2 \ge \mu_1 (g;k) \vol (T^2,g) \]

and equality holds if and only if the metric is flat. In particular $(k= \pm 1)$ we have 

\[ 4 \pi^2 \ge \mu_1 (g) \vol (T^2,g) \]

for the first positive eigenvalue of the Laplace operator in case that its eigenspace contains the representation $\Sigma (\pm 1)$.\\

\bigskip

Let us introduce the Hamiltonian operator $H_l$ defined by the Sturm-Liouville equation $(**)$ for $\lambda^2 =0$:

\[ H_l = - \frac{d^2}{dt^2} + 4 \pi l^2 - 4 \pi l \frac{h' (t)}{h(t)} - \frac{h(t)h''(t) - 2(h'(t))^2}{h^2(t)} . \]

\bigskip

{\bf Proposition 3:} {\it For $l \not= 0$, the Hamiltonian operators $H_l$  are strictly positive. $H_o$ is a non-negative operator.}\\

{\bf Proof:} A direct calculation yields the formula

\[ \int\limits^1_0 H_l \left( \frac{\varphi (t)}{h(t)} \right) \frac{\varphi (t)}{h(t)} dt = \int\limits^1_0 \left( 2 \pi l \frac{\varphi (t)}{h(t)} - \frac{\varphi' (t)}{h(t)} \right)^2 dt \ge 0 \]

where $\varphi (t)$ is any periodic function. The equation $2 \pi l \varphi (t) - \varphi'(t)=0$ does not admit a periodic, non-trivial solution in case $l \not= 0$. Consequently, $H_l$ is a strictly positive operator for $l \not= 0$.\\

{\bf Corollary:} {\it Fix an $S^1$-representation $\sum (l)$. Let $\lambda_1^2 (g,l)$ be the first eigenvalue of the Dirac operator on the two-dimensional torus $T^2$ with respect to the trivial spin structure and an $S^1$-invariant metric

\[ g=h^4(t) (dt^2 +dy^2) \]

such that the representation $\sum (l)$ occurs in the decomposition of the eigenspace. Then the multiplicity of $\sum (l)$  is one and the eigenspinor does not vanish anywhere $(l \not= 0)$.}\\

{\bf Proof:} Since $H_l$ is strictly positive, the eigenvalue $\lambda^2 (g)$ is the unique positive number $\lambda^2$ such that

\[ \mbox{inf} \, \, \mbox{spec} (H_l - \lambda^2 h^4)=0 . \]

\bigskip

The corresponding real solution of this Sturm-Liouville equation is unique and positive (see [8], page 207). \\

{\bf Corollary:} {\it For a fixed $S^1$-representation $\sum (l)$ denote by $\lambda_1^2 (g,l)$ the first eigenvalue of the Dirac operator such that the eigenspace $E(\lambda_1^2 (g,l))$ contains the representation $\sum (l)$. Then the inequality

\[ \lambda^2_1 (g,l) \le \frac{\D \int\limits^1_{0 \downsp} \frac{(2 \pi l \varphi (t) - \varphi' (t))^2}{\displaystyle h^2 (t)} \, dt}{\D \int\limits^{1^{\upsp}}_0 h^2 (t) \varphi^2 (t) dt}   \]

holds for any periodic function $\varphi (t)$.}\\

{\bf Proof:} Since $\mbox{inf} \, \, \mbox{spec} (H_l - \lambda_1^2 (g,l)h^4)=0$, we have

\[ \int\limits^1_0 H_l \left( \frac{\varphi}{h} \right) \frac{\varphi}{h} (t) dt  - \lambda_1^2 (g,l) \int\limits^1_0 h^2 (t) \varphi^2 (t) dt \ge 0 \]

for any periodic function $\varphi (t)$.\\

In case of the flat metric $g_o = dt^2+dy^2$ we have  $ \mu_1 (g_o) = \lambda_1^2 (g_o) = 4 \pi^2$ and

\[ E( \mu_1 (g_o))= \textstyle{\sum} (0) \oplus \sum (0) \oplus \sum (1) \oplus \sum  (-1)=E(\lambda_1^2 (g_o)) . \]

\bigskip

The spaces $\sum (\pm 1)$ correspond to the case that $k=l= \pm 1$ and are generated by the constant function. The two spaces $\sum (0)$ are generated by the functions $\sin (2 \pi t)$, $\cos ( 2 \pi t)$. \\

{\bf Notation:} We introduce now a few notations which will be used throughout this article. Let us consider a deformation $g_E = h^4_E (t) g_o$ of the flat metric $g_o$ depending on some parameter $E$. We assume that

\[ h^4_E (t) = h^4_E (1-t) \]

holds for all parameters of the deformation. The eigenvalues $\mu_1 (g_o)$ and $\lambda_1^2 (g_o)$ of multiplicity four split into three eigenvalues

\[ \mu_1 (g_o) \mapsto  \{ \mu_1 (E), \mu_2 (E), \mu_3 (E) \} \quad , \quad  \lambda_1^2 (g_o) \mapsto \{ \lambda_1^2 (E), \lambda_2^2 (E), \lambda_3^2 (E) \} . \]

The eigenvalue $\mu_3 (E)$ corresponds to the case that $k= \pm 1$, has multiplicity two, and its eigenfunction is a deformation of the constant function. The eigenvalues $\mu_1 (E) \not= \mu_2 (E)$ correspond to solutions of the Sturm-Liouville equation $(*)$ and their eigenfunctions are deformations of $\sin \, (2 \pi t)$ and $\cos \, (2 \pi t)$, respectively. The situation is different for the Dirac equation: there, according to Proposition 2,  the trivial $S^1$-representation $(l=0)$ yields one eigenvalue $\lambda_1^2 (E)$ of multiplicity two and the non-trivial representations $(l=\pm 1)$ define in general two distinct eigenvalues $\lambda_2^2 (E)$, $\lambda_3^2 (E)$ of multiplicity one. However, in case $h_E (t)=h_E (1-t)$, the spectral functions $\lambda_2^2 (E)$ and $\lambda_3^2 (E)$ coincide. Obviously, for small values $E \approx 0$ we have

\[ \mu_1 (g_E)= \min \, \{ \mu_1 (E), \mu_2 (E), \mu_3 (E) \} \quad , \quad  \lambda_1^2 (g_E)= \min \, \{ \lambda_1^2 (E), \lambda_2^3 (E), \lambda_3^2 (E) \} . \]

We will compute the first and second variation of $\mu_{\alpha} (E)$ and $\lambda_{\alpha}^2 (E)$ at $E=0$. For this purpose we introduce the following notation: Let $A$ be a function depending both on $E$ and $t$. Then $\dot{A}$ denotes the derivative with respect to $E$ and $A'$ the derivative with respect to $t$. Moreover, we expand the function $h^4_E (t)$ in the form

\[ h^4_E (t)= 1+ E H(t) + E^2 G(t) + {\cal O} (E^3) . \] \bigskip

\newcommand{\mpunkt}{\dot{\mu}}
\newcommand{\lpunkt}{\dot{\lambda^2}}

{\bf Theorem 1:} {\it Consider a deformation

\[ g_E = (1+EH(t)+E^2 G(t)+ {\cal O} (E^3)) {g_o} = h^4_E (t) g_o \]

of the flat metric on the torus $T^2$ such that $h_E^4 (t) = h_E^4 (1-t)$. Moreover, suppose that for $E \not= 0$ and $k=0$ the eigenvalues $\mu_1 (E) \not= \mu_2 (E)$ are simple eigenvalues of the Sturm-Liouville equation $(*)$. Then}

\begin{itemize}
\item[a.)] $\mpunkt_1  (0)= - 8 \pi^2 \D \int\limits_0^1 H(t) \sin^2 \, (2 \pi t) dt \, , \quad  \mpunkt_2 (0)= - 8 \pi^2 \D \int\limits^1_0 H(t) \cos^2 \, (2 \pi t) dt$\\
$\mpunkt_3 (0)= - 4 \pi^2 \D \int\limits^1_0 H(t) dt.$
\item[b.)] $\dot{\lambda}_1^2 (0)= \dot{\lambda}_2^2 (0)= \dot{\lambda}_3^2 (0)= - 4 \pi^2 \D \int\limits^1_0 H(t) dt$.
\end{itemize}

In particular, we obtain

\[ \mpunkt_1 (0) + \mpunkt_2 (0) = 2 \mpunkt_3 (0) = 2 \dot{\lambda}_{\alpha}^2 (0) . \] \bigskip

{\bf Corollary:} {\it Suppose that the deformation

\[ g_E =( 1+ EH(t)+ E^2 G(t) + {\cal O} (E^3)) g_o \]

of the flat metric $g_o$ satisfies the condition

\[ H(t)=H(1-t) \]

as well as 

\[ \D \int\limits^1_0 H(t) \sin^2 \, (2 \pi t)dt \not= \int\limits^1_0 H(t) \cos^2 (2 \pi t) dt . \]

Then, for all parameters $E \not= 0$ near zero we have the strict inequality}

\[ \mu_1 (g_E) < \lambda_1^2 (g_E) . \]

Next we compute the second variation of our spectral functions under the assumption that the first variation is trivial.\\

{\bf Theorem 2:} {\it Consider a deformation

\[ g_E =(1+EH(t) + E^2 G(t) + {\cal O} (E^3)) g_o = h^4_E (t) g_o \]

of the flat metric $g_o$ on the torus $T^2$ and suppose that the conditions

\[ h^4_E (t)=h^4_E (1-t) \]

and 

\[ \D \int\limits^1_0 H(t) \sin^2 \, (2 \pi t) dt = \int\limits^1_0 H(t) \cos^2 \, (2 \pi t) dt =0 \]

are satisfied. Moreover, suppose that for $E \not= 0$ and $k=0$ the eigenvalues $\mu_1 (E) \not= \mu_2 (E)$ are simple eigenvalues of the Sturm-Liouville equation $(*)$. Then

\newcommand{\doppel}{\ddot{\mu}}

\begin{itemize}
\item[a.)] $\doppel_1 (0)= - 16 \pi^2 \D \int\limits^1_0 G(t) \sin^2 \, (2 \pi t) dt - 16 \pi^2 \int\limits^1_0 H(t) C(t) \sin \, (2 \pi t) dt $,\\

where $C(t)$ is the periodic solution of the differential equation

\[ C'' (t) = - 4 \pi^2 H(t) \sin \, (2 \pi t) - 4 \pi^2 {C} (t) . \]

\item[b.)] $\doppel_2 (0) = - 16 \pi^2 \D \int\limits^1_0 G(t) \cos^2 (2 \pi t) - 16 \pi^2 \int\limits^1_0 H(t) {C} (t) \cos (2 \pi t) dt$,\\

 where ${C} (t)$ is the periodic solution of the differential equation

\[ C'' (t) = - 4 \pi^2 H(t) \cos (2 \pi t) - 4 \pi^2 {C} (t) . \]

\item[c.)] $\doppel_3 (0) = - 8 \pi^2 \D \int\limits^1_0 G(t) dt - 8 \pi^2 \int\limits^1_0 H(t) {C} (t)dt$,\\

where ${C} (t)$ is the periodic solution of the differential equation

\[ {C}'' (t) = - 4 \pi^2 H(t) . \]

\item[d.)] $\ddot{\lambda}_1^2 (0)=  - 8 \pi^2 \D \int\limits^1_0 G(t)dt + 2 \pi^2 \D \int\limits^1_0 H^2(t)dt$\\

\item[e.)]  $\ddot{\lambda}^2_2 (0)  =  \ddot{\lambda}^2_3 (0)  =  - 8 \pi^2 \D \int\limits^1_0 G(t) dt + 4 \pi^2 \int\limits^1_0 H^2 (t) dt - 8 \pi^2 \int\limits^1_0 H(t) C(t)dt  $\\
\mbox{} \hspace{2.5cm} $- 2 \pi \D \int\limits^1_0 H' (t) C(t) dt $,\\

where ${C} (t)$ is the periodic solution of the differential equation

\[ {C}''(t) = -  4 \pi^2 H(t)  -  \pi H' (t) . \]

\end{itemize}}

{\bf Proof of Theorem 1 and Theorem 2:} The formulas for the derivatives of $\lambda_1^2 (E)$ are a direct consequence of Proposition 2. We will prove the variation formulas for $\lambda^2_3$ and just remark that one can investigate the other spectral functions in a similar way. Moreover, since all the calculations we make are up to order two with respect to $E$,  we may assume for simplicity that

\[ h^4_E (t) = 1 + EH (t) + E^2 G(t) . \]

We compute

\[ \frac{h_E h_E'' - 2 (h'_E)^2}{h^2_E} = \frac{1}{4} E \frac{H'' + EG''}{(1+EH+E^2G)} - \frac{5}{16} E^2 \frac{(H' + EG')^2}{(1+EH+E^2G)^2} \]

and, consequently,  we obtain the formulas

\begin{eqnarray*}
 \frac{d}{dE} \left( \frac{h_E h_E'' - 2(h'_E)^2}{h^2_E} \right)_{E=0} &=& \frac{1}{4} H'' \quad , \quad \frac{d}{dE} \left( \frac{h'_E}{h_E} \right)_{E=0} = \frac{1}{4} H' \
\mbox{}\\
\mbox{}\\
\frac{d^2}{dE^2} \left( \frac{h_E h_E'' - 2(h'_E)^2}{h^2_E} \right)_{E=0} &=& \frac{1}{2} (G'' - H''H) - \frac{5}{8} (H')^2  . \end{eqnarray*}

\bigskip

The spectral function $\lambda_3^2 (E)$ is defined by a periodic solution $A_E(t)$ of the Sturm-Liouville equation

\[ A''_E (t) = - \lambda_3^2 (E)h^4_E (t) A_E(t) - \frac{h_E (t) h''_E (t) - 2(h'_E (t))^2}{h^2_E (t)} A_E (t) + 4 \pi^2 A_E (t) - 4 \pi \frac{h'_E (t)}{h_E (t)} A_E (t) \]

with the initial conditions $\lambda_3^2 (0) = 4 \pi^2, A_o (t) \equiv 1$. Therefore, we obtain 

\[ \dot{A}_o'' (t) = - \dot{\lambda}^2_3 (0) - 4 \pi^2 H(t) - 4 \pi^2 \dot{A}_o (t) - \frac{1}{4} H'' (t) + 4 \pi^2 \dot{A}_o (t) -  \pi H' (t) \]

in this case and, consequently,

\newcommand{\enull}{\int\limits^1_0}

\[ \dot{\lambda}_3^2 (0)= - 4 \pi^2 \enull H(t)dt . \]

\bigskip

Let us now compute the second variation in case that $\dot{\lambda}_1^2 (0)=0= \dot{\lambda}^2_3 (0)$. We differentiate the Sturm-Liouville equation twice at $E=0$: 

\begin{eqnarray*}
\ddot{A}_o'' (t) &=& - \ddot{\lambda}_3^2 (0) - 8 \pi^2 G(t) - 8 \pi^2 H(t) \dot{A}_o (t) + \frac{5}{8} (H' (t))^2  - \frac{1}{2} \Big(G'' (t) - H'' (t) H(t) \Big)\\
&&  - \frac{1}{2} H'' (t) \dot{A}_o (t) - 2 \pi H' (t) \dot{A}_o (t) - 4 \pi \frac{d}{dt} \left( \frac{d^2}{dE^2} (\mbox{ln} \, (h_E (t)))_{E=0} \right) .  
\end{eqnarray*}

\bigskip

Then we obtain

\begin{eqnarray*}
 \ddot{\lambda}^2_3 (0) &=& - 8 \pi^2 \enull G(t) dt - 8 \pi^2 \enull H(t) \dot{A}_o (t) dt  + \left( \frac{5}{8} - \frac{1}{2} \right) \enull (H' (t))^2 dt\\
&&  - \frac{1}{2} \enull H'' (t) \dot{A}_o (t) dt - 2 \pi \int\limits^1_0 H' (t) \dot{A}_o (t) dt   . 
\end{eqnarray*}

\bigskip

Since $\dot{A}_o(t)$ is a solution of the differential equation

\[ \dot{A}_o'' (t) = - 4 \pi^2 H(t) - \frac{1}{4} H'' (t)  -  \pi H' (t) ,  \]

we have

\begin{eqnarray*}
- \frac{1}{2} \enull H'' (t) \dot{A}_o (t) dt &=& - \frac{1}{2} \enull H(t) \dot{A}_o'' (t) dt =  + \frac{1}{2} \enull H(t) \left( 4 \pi^2 H(t) + \frac{1}{4} H'' (t) \right) dt \\
&=& 2 \pi^2 \enull H^2 (t) dt - \frac{1}{8} \enull (H' (t))^2 dt 
\end{eqnarray*}

and, consequently, we obtain

\begin{eqnarray*}
 \ddot{\lambda}^2_3 (0) &=& - 8 \pi^2 \enull G(t) dt + 2 \pi^2 \enull H^2 (t)dt - 8 \pi^2 \enull H(t) \dot{A}_o (t) dt - 2 \pi \int\limits^1_0 H' (t) \dot{A}_o (t) dt \\
&=&  - 8 \pi^2 \int\limits^1_0 G(t) dt + 4 \pi^2 \int\limits^1_0 H^2 (t) dt - 8 \pi^2 \int\limits^1_0 H(t) C(t)dt  - 2 \pi \int\limits^1_0 H' (t) C(t) dt  , 
\end{eqnarray*}

where ${C} (t) := \dot{A}_o (t) + \frac{1}{4} H(t)$ is a solution of the differential equation

\[ {C}''(t) = -  4 \pi^2 H(t)  -  \pi H' (t) . \]

\bigskip

{\bf Corollary:} {\it \, \, $\ddot{\lambda}_3^2 (0) = \ddot{\mu}_3 (0) + 2 \pi^2 \D \int\limits^1_0 H^2 (t) dt$.\\

In particular, for all parameters $E \not= 0$ near zero we have the inequality}

\[ \mu_3 (E) < \lambda_3^2 (E) . \]

{\it Moreover, the first positive eigenvalue $\mu_1 (g_E)$ of the Laplace operator is always smaller then the corresponding eigenvalue $\lambda_1^2 (g_E)$ of the Dirac operator for any metric $g_E$ near $E \approx 0$, i.e.,}

\[ \mu_1 (g_E) < \lambda_1^2 (g_E) . \]

\bigskip

{\bf Remark:}  The explicit formulas in { Theorem 1} and {Theorem 2} can be generalized to the case of an arbitrary conformal deformation. Fix a Riemannian metric $g_o$ on a surface $M^2$ and consider the deformation

\[ g_E = (1 + EH + E^2G + {\cal O} (E^3)) g_o \]

of the metric. Moreover, suppose that $\mu_1 (E)$ is the deformation of the eigenvalue of the Laplace operator and $f_E$ is the corresponding family of eigenfunctions. Then the following  formulas hold:

\begin{itemize}
\item[a.)] $\dot{\mu}_1 (0)= - \mu_1 (0) \, \, \frac{\D \int\limits_{M^2 \downsp} H f_o^2 dM^2_o}{\D \int\limits_{M^2}^{\upsp} f_o^2 dM^2_o}$
\item[b.)] $\ddot{\mu}_1 (0)= - 2 \mu_1 (0) \, \, \frac{ \D \int\limits_{M^2 \downsp} (G f_o^2 + C f_o)dM^2_o}{\D \int\limits_{M^2}^{\upsp}  f_o^2 dM^2_o}$ , 
\end{itemize}

where the function $C$ is the solution of the differential equation

\[ \Delta_o C= \mu_1 (0) H f_o + \mu_1 (0) C . \]

The corresponding expression for the variation of the eigenvalue $\lambda_1 (E)$ of the Dirac operator can also be computed:
\medskip

\begin{itemize}
\item[c.)] $\dot{\lambda}_1 (0)= - \lambda_1 (0) \, \, \frac{\D  \int\limits_{M^2 \downsp} H \cdot |\psi_o|^2 dM^2_o}{ \D 2 \int\limits_{M^2}^{\upsp} |\psi_o|^2 dM^2_o}$
\end{itemize}

and a similar formula holds  for the second variation.\\

Once again, a similar, though even more intricate computation yields the fourth variation of $\lambda_3^2 (E)$ under the assumption that all previous variations of $\lambda_3^2 (E)$ vanish. This is needed for the discussion of the example in Section 5.\\

{\bf Theorem 3:}  {\it Consider a deformation

\[ g_E =(1+EH(t)) g_o \]

of the flat metric $g_o$ on the torus $T^2$ and suppose that the following conditions are satisfied:

\begin{itemize}
\item[a.)] $H(t) =H (1-t);$
\item[b.)] $\dot{\lambda}_3^2 (0)= \ddot{\lambda}_3^2 (0) = \stackrel{\boldmath \cdots \unboldmath}{\lambda} \! \!\mbox{}^2_3 \,  (0)=0$ . 
\end{itemize}

Then the fourth derivative $[\lambda_3^2 (0)]^{\mbox{\tiny (IV)}}$ of the spectral function $\lambda_3^2 (E)$ at $E=0$ is given by the formula

\begin{eqnarray*}
[ \lambda_3^2 (0)]^{\mbox{\tiny (IV)}} &=& 6 \int\limits^1_0 H^3 (t) H'' (t)dt + \frac{45}{2} \int\limits^1_0 H^2 (t) (H' (t))^2 dt  - \int\limits^1_0 \Big(16 \pi^2 H(t) + H''(t) \Big) C_3 (t) dt \\
&& + \int\limits^1_0 \left( \frac{5}{2} (H' (t))^2 + 2H(t)H''(t) \right) C_2 (t)dt  \\
&& - \int\limits^1_0 \Big(15 (H'(t))^2 + 6 H''(t)H(t) \Big) H(t) C_1 (t) dt , 
\end{eqnarray*}

where the functions $C_1 (t), C_2(t), C_3 (t)$ are periodic solutions of the equations

\begin{eqnarray*}
C_1'' (t) &=& - 4 \pi^2 H(t) - \frac{1}{4} H'' (t) - \pi H' (t) \\
C_2'' (t) &=&  \frac{1}{2} H(t) H'' (t) + \frac{5}{8} (H' (t))^2 + 2 \pi H' (t) H(t)  -  \left( 8 \pi^2 H(t) + \frac{1}{2} H''(t)+ 2 \pi H' (t) \right) C_1 (t)\\
C_3'' (t) &=& - \left( \frac{3}{2} H'' (t) H^2 (t) + \frac{15}{4} H(t) (H' (t))^2 + 6 \pi H' (t) H^2 (t) \right)\\
&& + \left( \frac{3}{2} H(t) H''(t) + \frac{5}{8} (H' (t))^2 + 6 H(t) H' (t) \right) C_1 (t)  \\
&& - \left( 3 \pi H' (t) + \frac{3}{4} H'' (t) + 4 \pi^2 H(t) \right) C_2 (t) . 
\end{eqnarray*} }

{\bf Remark:} In the special case of $H'' (t)= - 16 \pi^2 H(t)$ the derivative $[\lambda_3^2 (0)]^{\mbox{\tiny (IV)}}$ does not depend on $C_3 (t)$ and the formulas become much simpler. Such a metric will be the object of Section 5.\\

\section{Examples}

\subsection{The variation \boldmath $g_E=(1+E \cos (2 \pi t))g_o$ \unboldmath}

The volume ${\vol} (T^2,g_E)=1$ of this variation of the flat metric $g_o$ is constant and all first derivatives at $E=0$ vanish since

\[ \int\limits^1_0 \cos (2 \pi t) \cos^2 (2 \pi t) = \int\limits^1_0 \cos (2 \pi t) \sin^2 (2 \pi t)=0 . \]

A computation of the second derivatives yields the following numerical values:

\[ \ddot{\mu}_1 (0) = - \frac{2}{3} \pi^2 \quad , \quad \ddot{\mu}_2 (0) = \frac{10}{3} \pi^2 \quad , \quad \ddot{\mu}_3 = - 4 \pi^2 \]

\[ \ddot{\lambda}_1^2 (0)= \pi^2 \quad , \quad \ddot{\lambda}_2^2 (0)  = \ddot{\lambda}_3^2 (0) = - 3 \pi^2 . \]

In particular, we obtain

\[ \mu_1 (g_E) < \lambda_1^2 (g_E) \]

for all parameters $E \not= 0$ near zero. The eigenspinor corresponding to the minimal positive eigenvalue of the Dirac operator does not vanish anywhere (Figure 1).\\



\begin{center}

\[
\epsfig{figure=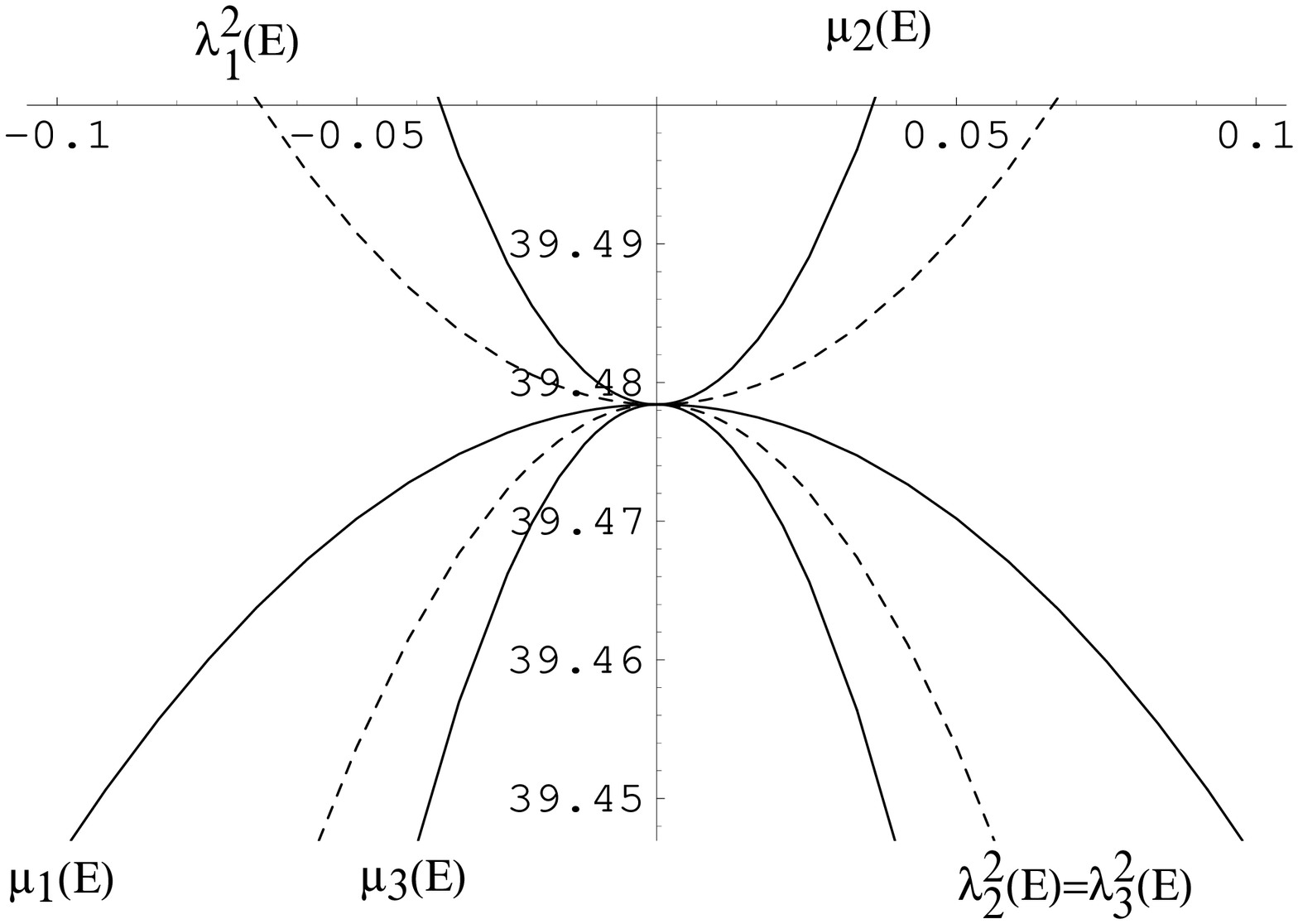,width=10cm}
\]

(Figure 1)
\end{center}

\vspace{1cm}

\subsection{The Mathieu deformation \boldmath $ g_E =( 1 + E \cos (4 \pi t)) g_o$ \unboldmath of the flat metric}

This deformation of the flat metric again preserves the volume, and the Laplace equation essentially reduces to the classical Mathieu equation $u'' (x) + (a + 16 q \cos (2x)) u(x)=0$. In this case the first variation is trivial only for the Dirac equation. Indeed, we have

\[ \dot{\mu}_1 (0) = 2 \pi^2  \quad , \quad \dot{\mu}_2 (0) = - 2 \pi^2 \quad , \quad \dot{\mu}_3 (0)=0 \]

\[ \dot{\lambda}_1^2 (0) = \dot{\lambda}^2_2 (0) = \dot{\lambda}_3^2 (0)= 0 . \]

Even for the Mathieu deformation we  conclude that 

\[ \mu_1 (g_E) < \lambda^2_1 (g_E) \]

for $E \not=0$ near zero. A computation of the second derivatives yields the following result (Figure 2):

\[ \ddot{\mu}_3 (0) = - \pi^2 \quad , \quad  \ddot{\lambda}_1^2 (0) = \pi^2 \quad , \quad  \ddot{\lambda}^2_2 (0) = \ddot{\lambda}_3^2 (0)= 0 . \]

For a detailed discussion of this metric, we refer to the next section.\\

\bigskip

\begin{center}

\[
 \epsfig{figure=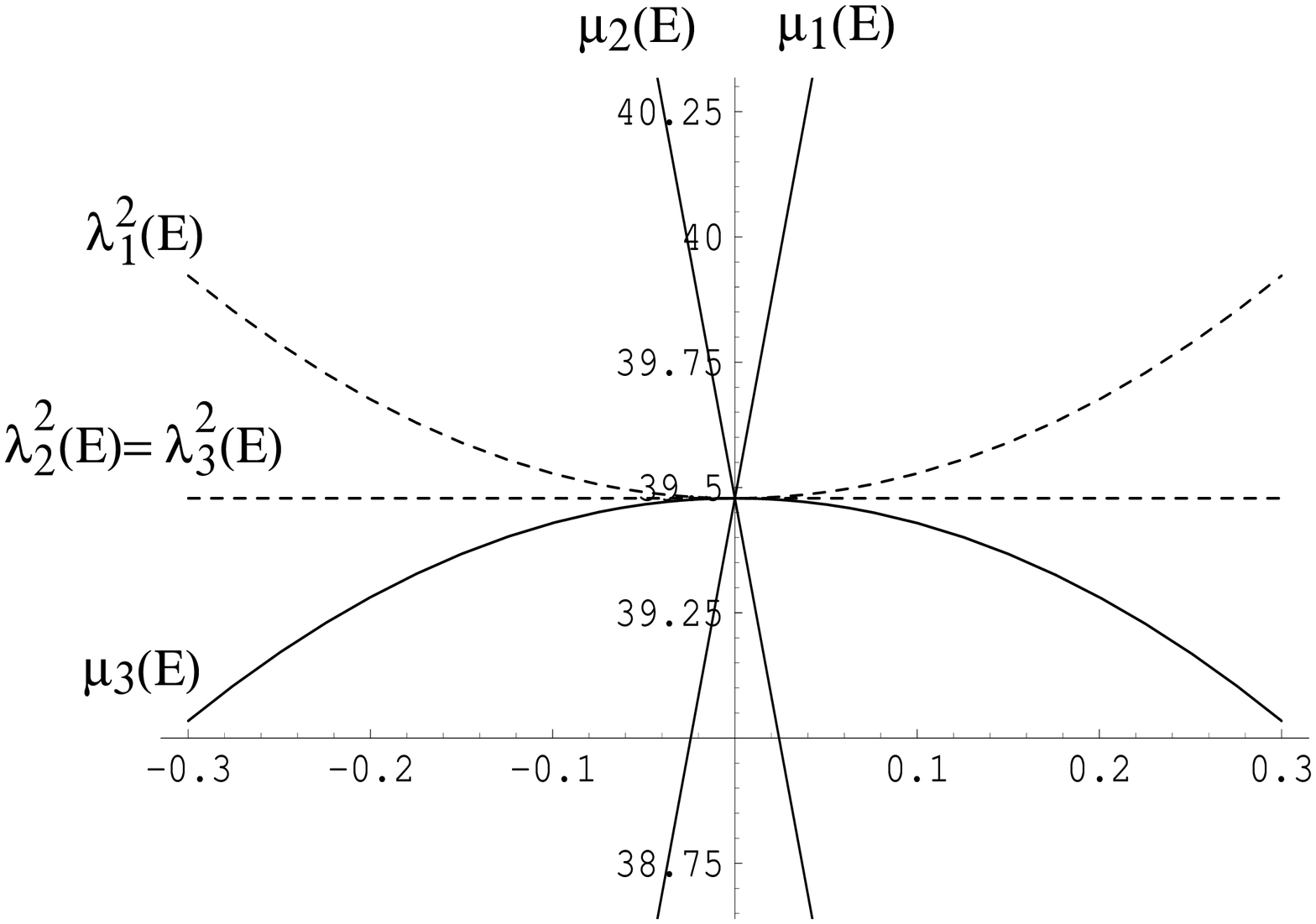,width=10cm} \hspace{1cm} \mbox{}
\]

(Figure 2)
\end{center}

\vspace{1cm}

\subsection{The variation \boldmath $g_E =( 1 + E \cos (2 \pi Nt)) g_o, \, \, \, N \ge 3$ \unboldmath}

Since

\[  \int\limits^1_0 \cos (2 \pi Nt) \cos^2 ( 2 \pi t) dt= \int\limits^1_0  \cos (2 \pi Nt) \sin^2 (2 \pi t) dt=0 \]

for $N \ge 3$, the first variations of our spectral function vanish. We compute the second variation using the algorithm in { Theorem 2}:

\[ \ddot{\mu}_1 (0) = \ddot{\mu}_2 (0) = - \frac{4 \pi^2}{N^2 - 4} \quad , \quad \ddot{\mu}_3 (0) = - \frac{4 \pi^2}{N^2} \]

\[ \ddot{\lambda}_1^2 (0) =  \pi^2 \quad , \quad  \ddot{\lambda}^2_2 (E)= \ddot{\lambda}_3^2 (0)= \left( 1- \frac{4}{N^2} \right) \pi^2 . \]

\bigskip

In particular, we obtain again

\[ \lambda_1^2 (g_E) > \mu_1 (g_E) \]

for all parameters $E \not= 0$ near zero (Figure 3). \\

\bigskip

\begin{center}

\[
\epsfig{figure=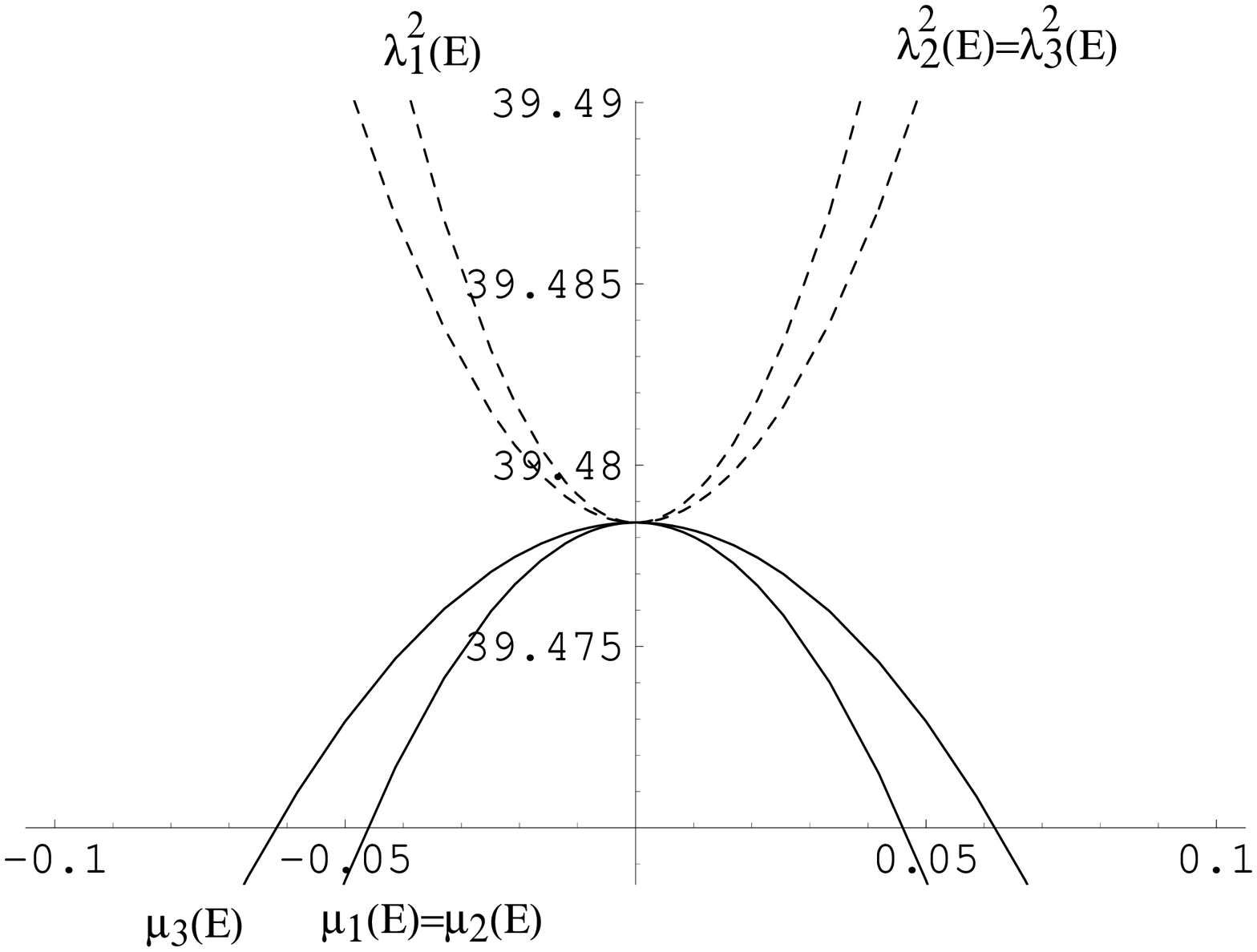,width=10cm}
\]

(Figure 3)
\end{center}

\vspace{1cm}

\section{The Mathieu deformation of the flat metric}

In the previous examples, the deformation

\[ g_E = (1+E \cos (4 \pi t)) g_o \]

of the flat metric $g_o$ plays  an exceptional role, because the derivatives $\dot{\mu}_1 (0), \dot{\mu}_2 (0) \not= 0$ are non-zero. Therefore, we study  the behaviour of the first positive eigenvalue for the Laplace and Dirac operator in more detail. First of all, the lower bound

\[ \frac{4 \pi^2}{h^4_{\max}} \le \mu_1 (E) , \lambda_1^2 (E) \]

yields the estimate

\[ \frac{4 \pi^2}{1 + |E|} \le \mu_1 (E), \lambda_1^2 (E) \]

for all parameters $-1 <E<1$.  In case of the function $f(t) = \sin (2 \pi t)$ the upper bound $B^u_L (g,f)$ of Section 2 leads to the estimate

\[ \mu_1 (E) \le \frac{8 \pi^2}{2+|E|} , \]

i.e., for all parameters $-1<E<1$ the inequality

\[ \frac{4 \pi^2}{1 + |E|} \le \mu_1 (E) \le \frac{8 \pi^2}{2 + |E|} \]

holds. On the other hand, for the Dirac operator the function

\[ f(t) = (1+E \cos (4 \pi t))^{\frac{1}{4}} \sin (2 \pi t) \]

gives an upper bound $B_D^u (g_E,f)$ for its first eigenvalue with the property

\[ \lim\limits_{E \to - 1} B^u_D (g_E, f) = 5 \pi^2 .  \]

\bigskip

We will thus investigate the limits $\lim\limits_{E \to -1} \mu_1 (E)$ as well as $\lim\limits_{E \to - 1} \lambda_1^2 (E)$. The eigenvalue $\mu_1 (E)$ is related with a periodic solution of the Sturm-Liouville equation

\[ A'' (t) = - \mu_1 (E) \Big(1+E \cos ( 4 \pi t) \Big) A(t) + 4 \pi^2 k^2 A (t)  , \]

where $k=0, \pm 1$ (see Proposition 1). Let us introduce the function $B(x) := A \left( \frac{1}{2 \pi} x \right)$ where $0 \le x \le 2 \pi$. Then the Sturm-Liouville equation is equivalent to the classical Mathieu equation

\[ B'' (x) + (a + 16q \cos (2x)) B(x) =0 , \]

where the parameters $a$ and $q$ are given by

\[ a= \frac{\mu_1 (E)}{4 \pi^2} - k^2 \quad , \quad q= \frac{E \mu_1 (e)}{16(4 \pi^2)} \quad , \quad k=0, \pm 1 . \]

\bigskip

For $E \to - 1$ the parameters of the Mathieu equation are related by

\[ a= - 16 q - k^2 \quad , \quad k=0, \pm 1 . \]

\bigskip

Using the estimates for $\mu_1 (E)$ we obtain 

\[ 2 \pi^2 \le \lim\limits_{E \to -1} \mu_1 (E) \le \frac{8}{3} \pi^2 , \]

\[ \mbox{i.e.,} \, \, - \frac{1}{24} \le q \le - \frac{1}{32} \quad  \mbox{in case} \, \, E= - 1 . \]

A numerical computation shows that, under these restrictions, the Mathieu equation has a unique periodic solution for $k=0$ and $q \approx 0,04113$. This solution $B(x)$ is the first Mathieu function $se_1 (x,q)$, which is the deformation of the function $\sin (x)$. Consequently, we have

\[ \lim\limits_{E \to -1} \mu_1 (E) = - 16 \cdot q \cdot 4 \pi^2 \approx 2,6323 \pi^2 . \]

\bigskip

The limits of the spectral functions $\mu_2 (E)$ and $\mu_3 (E)$ can be computed in a similar way:

\[ \lim\limits_{E \to -1} \mu_2 (E) \approx 1,79 \cdot (4 \pi^2) \quad , \quad \lim\limits_{E \to - 1} \mu_3 (E) \approx 0,9 \cdot (4 \pi^2) . \]

\bigskip

These limits correspond to the Mathieu functions $ce_1 (x,q)$ and $ce_o (x,q)$ for the parameters $q \approx - 0,112$ in case of $\mu_2 (E)$ and $q \approx - 0,056296$ in case of $\mu_3 (E)$. \\


\begin{center}

\[
\epsfig{figure=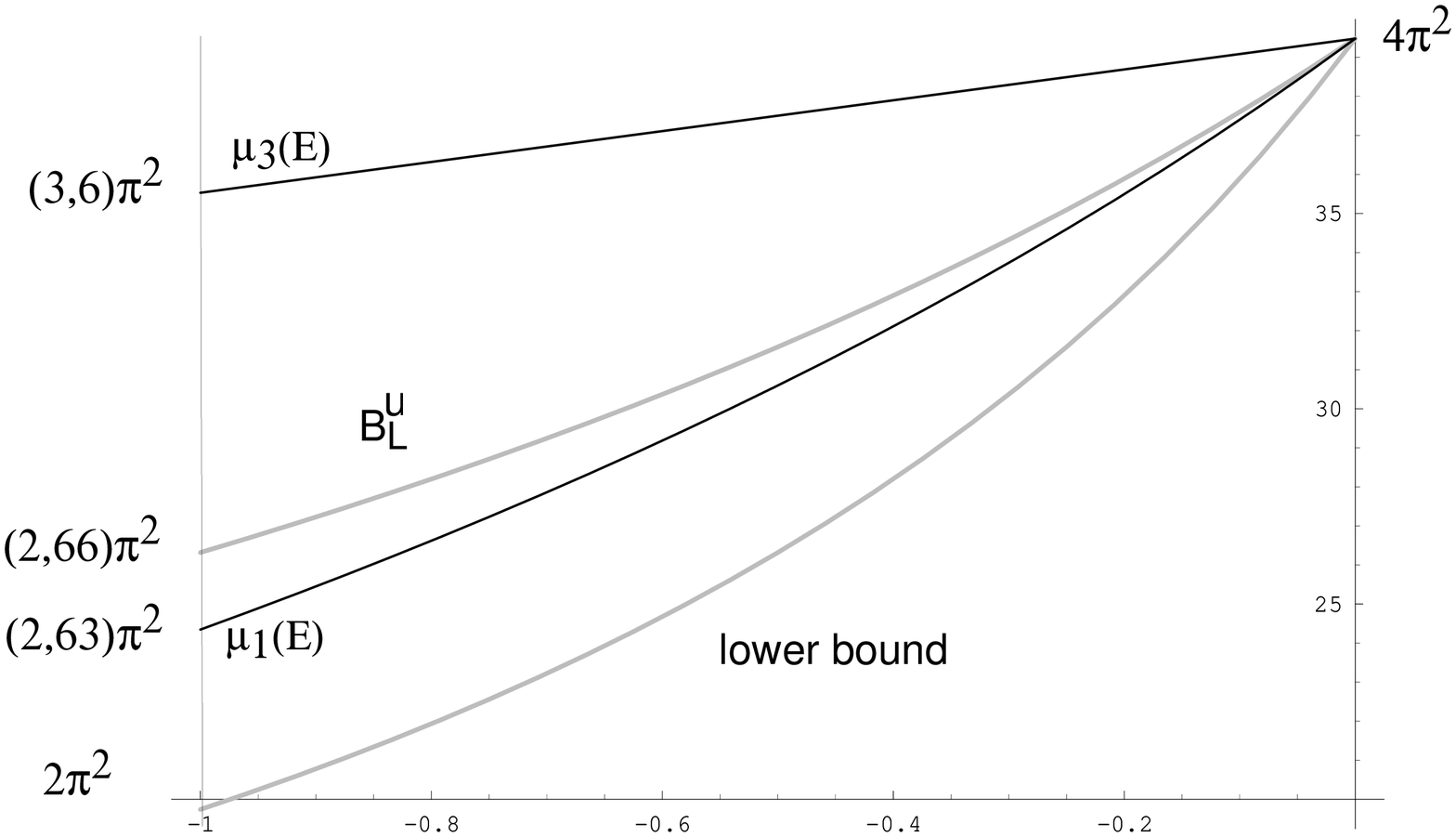,width=10cm}
\]

(Figure 4)
\end{center}

\vspace{1cm}


{\bf Approximation of the periodic solution for \boldmath $\mu_1 (E), E \to -1$ \unboldmath :}\\
\begin{center}
\begin{quote}
\begin{verbatim}
NDSolve[{y''[x] + 32(0.04113)(Sin[x])^2 y[x] == 0,
y[x] == 0 , y'[0] == 1} , y , {x , 0 , 10 Pi}]

Plot[Evaluate[y[x]/.% , {x , 0 , 10 Pi}]
\end{verbatim}
\end{quote}

\[
\epsfig{figure=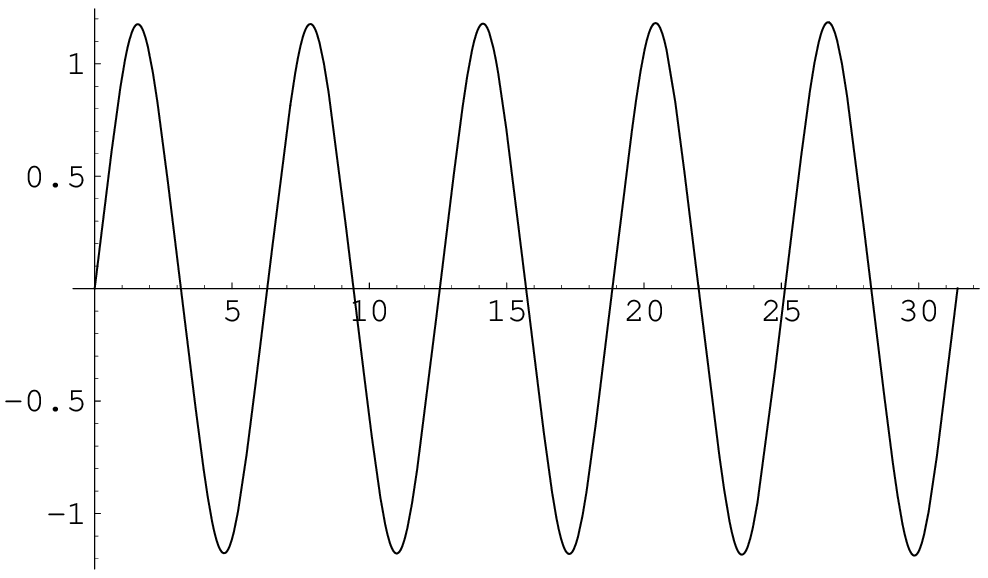,width=10cm}
\]

(Figure 5)
\end{center}

\vspace{1cm}

{\bf Approximation of the periodic solution for \boldmath $\mu_2 (E), E \to -1$ \unboldmath }:\\
\begin{center}
\begin{quote}
\begin{verbatim}
NDSolve[{y''[x] + 32(0.1112)(Sin[x])^2 y[x] == 0,
y[x] == 1 , y'[0] == 0} , y , {x , 0 , 10 Pi}]

Plot[Evaluate[y[x]/.% , {x , 0 , 10 Pi}]
\end{verbatim}
\end{quote}

\[
\epsfig{figure=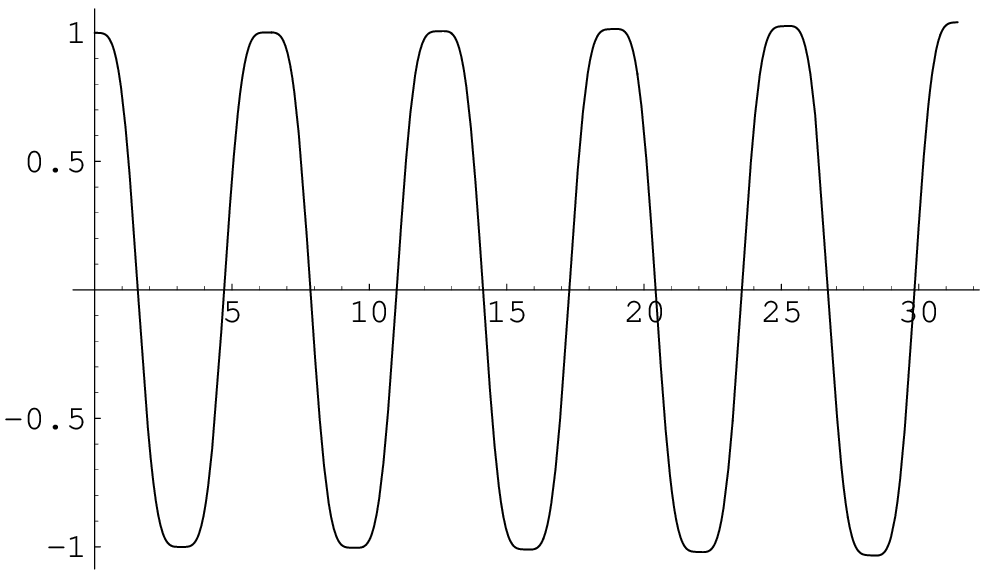,width=10cm}
\]

(Figure 6)
\end{center}

\vspace{1cm}

{\bf Approximation of the periodic solution for \boldmath $\mu_3 (E), E  \to -1$ \unboldmath :}\\
\begin{center}
\begin{quote}
\begin{verbatim}
NDSolve[{y''[x] + 32(0.056296)(Sin[x])^2 - 1) y[x] == 0,
y[0] == 1 , y'[0] == 0} , y , {x , 0 , 10 Pi}]

Plot[Evaluate[y[x]/.% , {x , 0 , 10 Pi}]
\end{verbatim}
\end{quote}

\[
\epsfig{figure=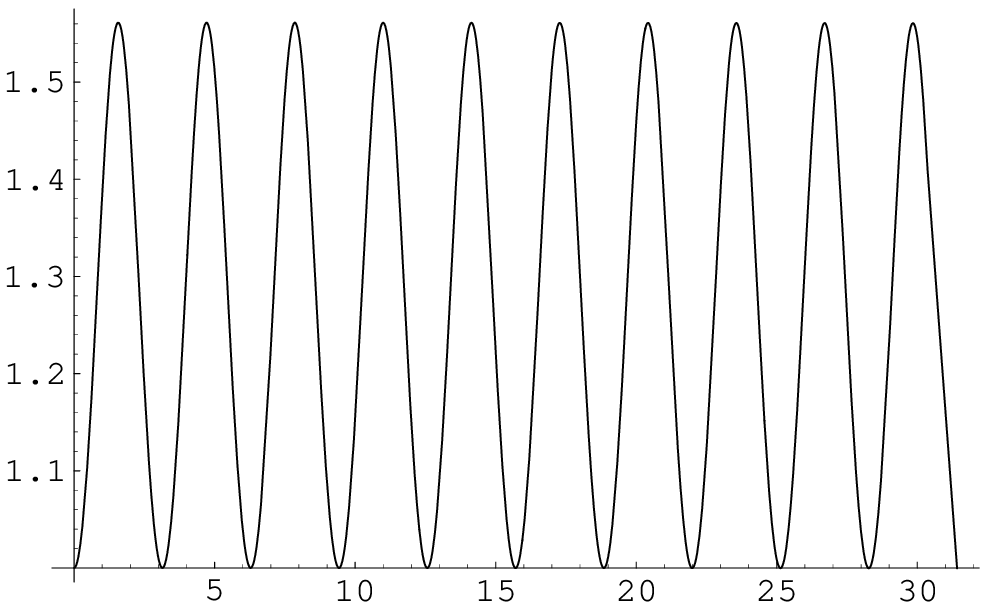,width=10cm}
\]

(Figure 7)
\end{center}

\vspace{1cm}

The eigenvalues $\lambda_{\alpha}^2 (E)$ of the Dirac operator are related with the periodic solutions of the Sturm-Liouville equation

\[ - A'' (t) = \left\{ \lambda^2 h^4 (t) + \frac{h(t) h'' (t) - 2 (h'(t))^2}{h^2(t)} - 4 \pi^2 l^2 \, + 4 \pi l \frac{h'(t)}{h(t)} \right\} A(t) .\]

\bigskip

For the Mathieu deformation we have

\[ \frac{h(t) h''(t) - 2(h'(t))^2}{h^2(t)} = - 4 \pi^2 E \frac{E + \cos (4 \pi t) + \frac{1}{4} E \sin^2 (4 \pi t)}{(1+E \cos (4 \pi t))^2} . \]

\bigskip

First we discuss the case that $l=0$. Then the first positive eigenvalue of the Dirac equation is given by

\[ \lambda^2 = \frac{4 \pi^2}{\left( \D \int\limits^1_0 h^2 (t) dt \right)^2} . \]

\bigskip

In case of the Mathieu deformation we obtain

\[ \lim\limits_{E \to -1} \int\limits^1_0 h^2 (t) dt = \int\limits^1_0 \sqrt{1- \cos (4 \pi t)} dt = \frac{2 \sqrt{2}}{\pi} \]

and, finally,

\[ \lim\limits_{E \to -1} \lambda^2 (E) = \frac{1}{2} \pi^4 \approx (4,92) \pi^2\] 

\bigskip

\begin{center}
\[
\epsfig{figure=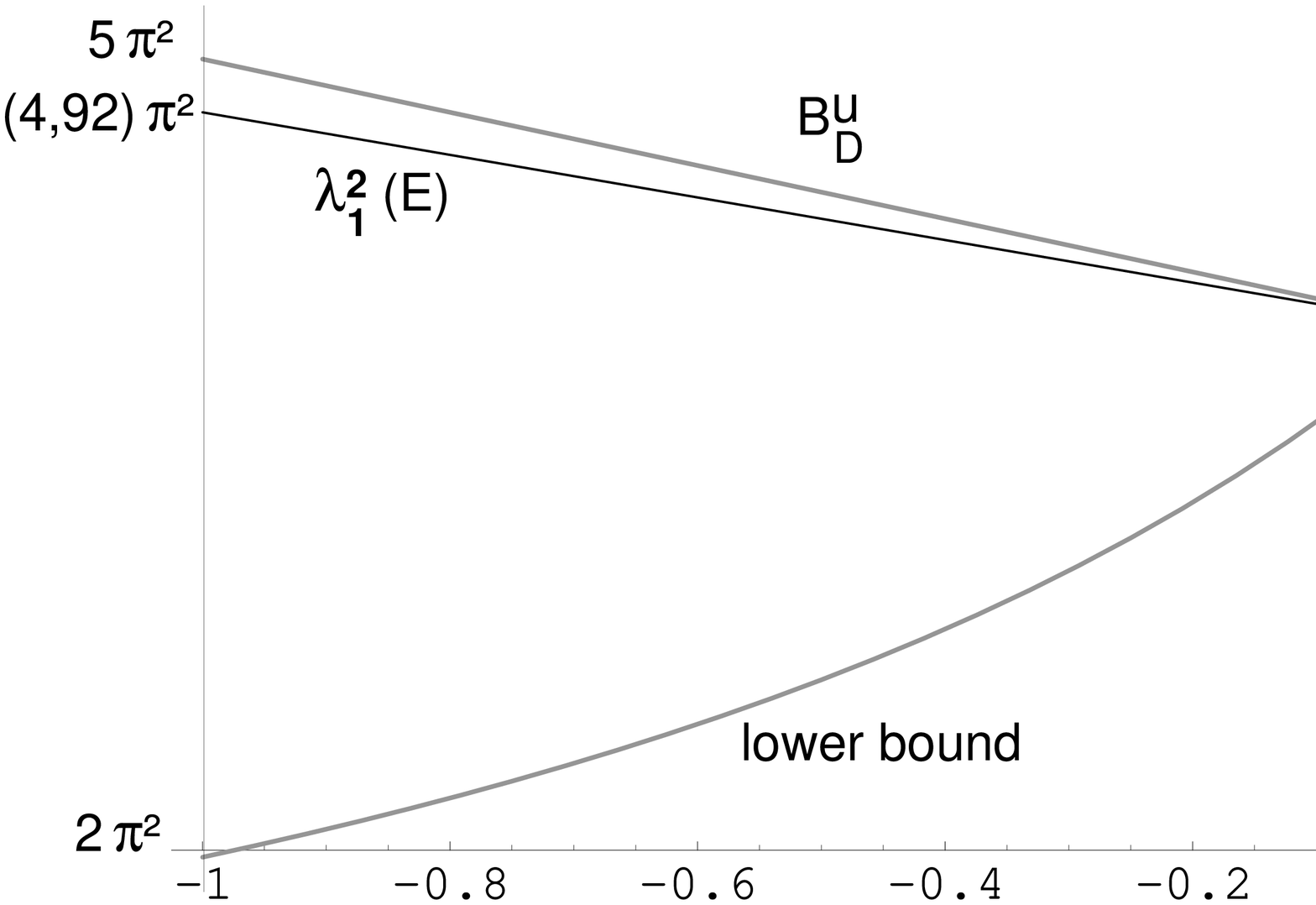,width=10cm}
\]

(Figure 8)
\end{center}

\vspace{1cm}

We now investigate the case $l=1$. Let us consider the Hamiltonian operator $H_E$ given by the Sturm-Liouville equation for $\lambda^2 =0$:

\begin{eqnarray*}
H_E = - \frac{d^2}{dt^2} + 4 \pi^2 - 4 \pi \frac{h'(t)}{h(t)} - \frac{h(t)h''(t) - 2(h' (t))^2}{h^2(t)} = - \frac{d^2}{dt^2} + p_E(t)  , 
\end{eqnarray*}

where the potential $p_E(t)$ is given by the formula

\[ p_E (t) = 4 \pi^2 + E \pi^2  \frac{4 \cos (4 \pi t) + 4 \sin (4 \pi t) + E \sin^2 (4 \pi t) + 2E (2+ \sin (8 \pi t))}{(1+E \cos ( 4 \pi t))^2} . \]

\bigskip

For all parameters $-1 < E \le 0$ the Hamiltonian operator $H_E$ is strictly positive (see Proposition 3). Consequently, the eigenvalue $\lambda^2_3 (E)$ is the first number such that

\[ \mbox{inf} \, \mbox{spec} \, (H_E - \lambda^2 (1+ E \cos (4 \pi t)))= 0 , \]

and the corresponding solution of the Sturm-Liouville equation

\[ A_E'' (t) = ( p_E(t) - \lambda_3^2 (E) (1+E \cos (4 \pi t)) ) A_E(t) \]

is unique and everywhere positive. In particular, the solution satisfies the condition

\[ A_E(t + \textstyle{\frac{1}{2}})=A_E (t) . \]

\bigskip

Since $A_E(t)$ is a positive periodic solution of the Sturm-Liouville equation, we obtain the condition

\[ \int\limits^1_0 (p_E(t) - \lambda_3^2 (E) (1+E \cos (4 \pi t))) dt >0 \]

and thus an upper bound for $\lambda_3^2 (E)$: 

\[  \lambda_3^2 (E)< \int\limits^1_0 p_E(t) dt . \]


\begin{center}
\[
\epsfig{figure=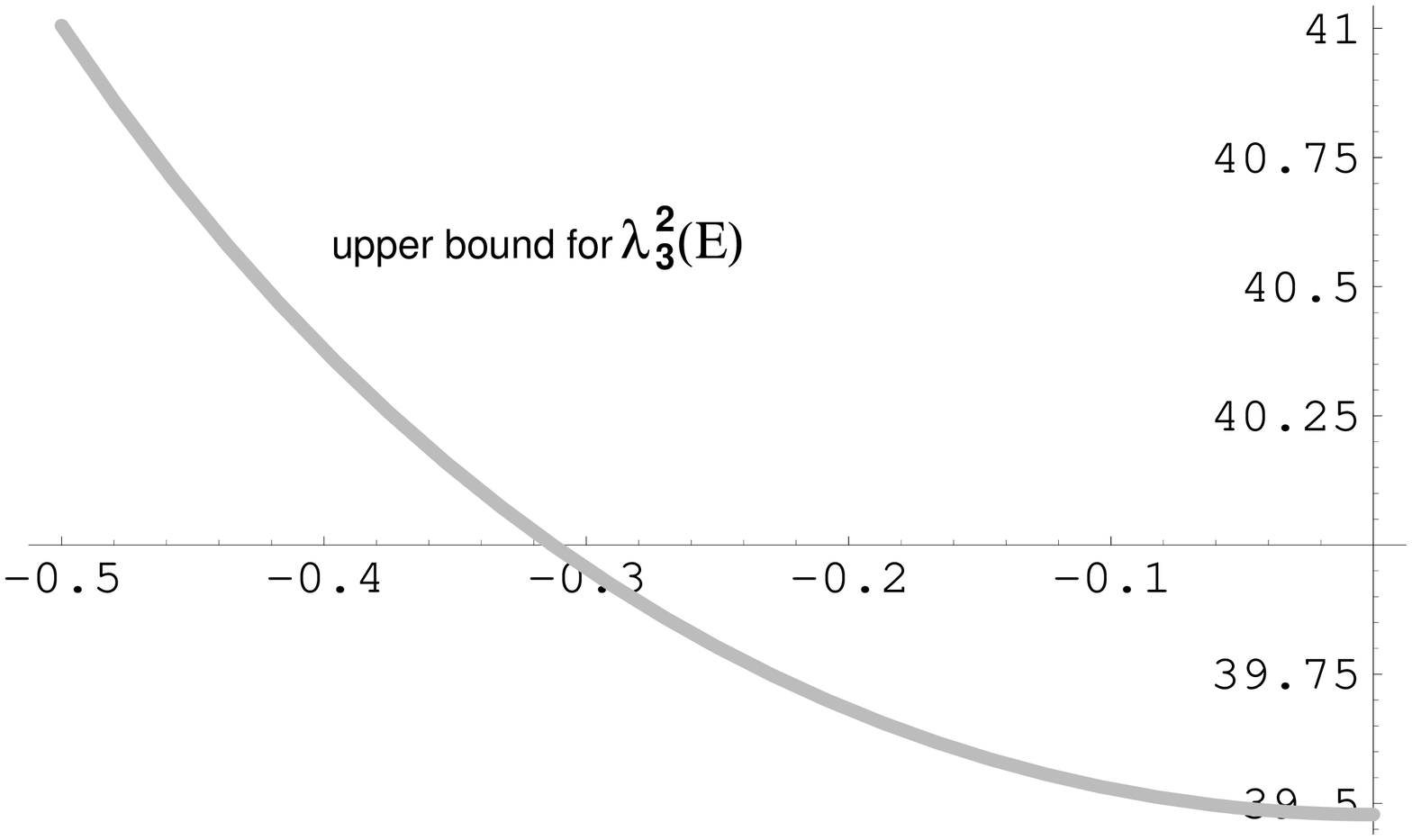,width=10cm}
\]

(Figure 9)
\end{center}

\vspace{1cm}

We notice that this upper bound for $\lambda_3^2 (E)$ grows  and  reflects, indeed,  the real behaviour of $\lambda_3^2 (E)$ near $E=0$. To see this, we use Theorem 3 to compute the fourth variation of this spectral function (the third variation vanishes since $\lambda_3^2 (E)$ has to be a symmetric function in $E$). One obtains the following result:

\[ [\lambda_3^2 (0)]^{\mbox{\tiny (IV)}} = \frac{27}{4} \pi^2 > 0 . \]

\bigskip

On the other hand, using well-known approximation techniques for Sturm-Liouville equations with periodic coefficients (see [9]) we can approximate $\lambda_3^2 (E)$ for a fixed parameter $E$. Indeed, one replaces the potential in the Sturm-Liouville equation by the first terms of its Fourier series. This reduces the computation of the approximative eigenvalue to a finite-dimensional eigenvalue problem. For example,  in case of $E=-0.3$ the mentioned methods yields the result

\[ \lambda_3^2 (-0.3) \approx 39.6733 . \]

\bigskip

Let us study the behaviour of the spectral function $\lambda_3^2 (E)$ for $E \to - 1$. More generally, denote by $\lambda^2 (E,l)$ the first eigenvalue of the Dirac operator such that the corresponding eigenspace contains an $S^1$-representation of weight $l$. In particular, we have $\lambda_3^2 (E) = \lambda^2 (E,1)= \lambda^2 (E,-1)$. We apply the Corollary of Proposition 3 to the function $h_E (t) = \sqrt[4]{1+E \cos ( 4 \pi t)}$ and conclude that

\[ \int\limits^1_0 \frac{(2 \pi l \varphi (t) - \varphi' (t))^2}{\sqrt{1+E \cos (4 \pi t)}} - \lambda^2 (E,l) \int\limits^1_0 \sqrt{1+E \cos (4 \pi t)} \varphi^2 (t) dt \ge 0 \]

holds for any periodic function $\varphi (t)$. Fix a test function $\varphi (t)$ and consider the limit $E \to -1$. Then we obtain the inequality

\[ \overline{\lim\limits_{E \to -1}} \lambda^2 (E,l) \le \frac{1}{2} \, \, \frac{\D \int\limits^1_{0} \frac{(2 \pi l \varphi (t) - \varphi' (t))^2}{|\sin (2 \pi t)|} dt}{\D \int\limits^1_0 |\sin (2 \pi t)| \varphi^2 (t) dt} . \]

We apply this estimate  to the function

\[ \varphi_l (t) = \frac{\cos (2 \pi t) + l \sin (2 \pi t)}{2(l^2 +1) \pi} . \]

\bigskip

Then $2 \pi l \varphi_l (t) - \varphi_l' (t) = \sin (2 \pi t)$ and we obtain the following\\

{\bf Proposition 4:} 

\[ \overline{\lim\limits_{E \to -1}} \lambda^2 (E,l) \le 6 \pi^2 \frac{(l^2+1)^2}{(1+2l^2)} . \]

\bigskip

{\bf Remark:} At $E=0$ we have $\lambda^2 (0,l)= 4 \pi^2 l^2$. On the other hand, for $l \ge 3$ the inequality

\[ 6 \cdot \frac{(l^2 +1)^2}{1+2l^2} < 4l^2 \]

holds, i.e.,\\

\mbox{} \hspace{4.7cm} $ \D  \overline{\lim\limits_{E \to -1}} \lambda^2 (E,l) < \lambda^2 (0,l)$ \hspace{2cm} $l \ge 3 . $\\

The latter inequality means that the eigenvalue $\lambda^2 (E,l)$ decreases for $E \to -1$ \, $(l \ge 3)$.\\


The behaviour of $\lambda_3^2 (E)= \lambda^2 (E,l)$ for $l=1$ is completely different. This spectral function increases for $E \to - 1$. Using the formula

\[ \lambda_3^2 (E)= \inf\limits_{\varphi > 0} \frac{\D \int\limits^1_0 \frac{\D (2 \pi l \varphi (t) - \varphi' (t))^2}{\D \sqrt{1+E \cos (4 \pi t)}}}{\D \int\limits^1_0 \sqrt{1+E \cos (4 \pi t)} \varphi^2 }\]

we can approximate the positive minimizing  {\it Mathieu spinor} $\mbox{MS} \, (E,t)$ of topological index $l=1$   by expanding it in its Fourier series. We thus obtain for example:\\

\bigskip

{\boldmath $E= - 0.9$\unboldmath}:  \quad $\lambda_3^2 (-0.9) \approx 40.1464$\\

\begin{verbatim} 
MS(-0.9,t)=(Sqrt[Sqrt[1+ (-0.9)Cos[4 Pi t]]]) Sqrt[1+ (0.44)Sin[4 Pi t] 
    + (0.15)Cos[4 Pi t] + (0.09)Sin[8 Pi t] + (0.17)Cos[8 Pi t]
    + (0.028)Sin[16 Pi t] + (0.051)Cos[16 Pi t] + (0.051)Sin[12 Pi t] 
    + (0.085)Cos[12 Pi t] + (0.016)Sin[20 Pi t] + (0.026)Cos[20 Pi t]
    + (0.01)Sin[24 Pi t] + (0.014)Cos[24 Pi t] + (0.005)Sin[28 Pi t] 
    + (0.007)Cos[28 Pi t] + (0.0033)Sin[32 Pi t] + (0.0044)Cos[32 Pi t]] \end{verbatim} 

\bigskip

{\boldmath $E= - 0.95$\unboldmath}:  \quad $\lambda_3^2 (-0.9) \approx 44.6024$\\

\begin{verbatim} 
MS(-0.95,t)=(Sqrt[Sqrt[1+ (-0.95)Cos[4 Pi t]]])Sqrt[1+ (0.585)Sin[4 Pi t] 
    + (0.049)Cos[4 Pi t] + ( 0.1)Sin[8 Pi t] + (0.08)Cos[8 Pi t] 
    + (0.063)Sin[12 Pi t] + (0.063)Cos[12 Pi t] + (0.041)Sin[16 Pi t] 
    + (0.04)Cos[16 Pi t]  + (0.026)Sin[20 Pi t] + (0.026)Cos[20 Pi t]  
    + (0.017)Sin[24 Pi t] + (0.018)Cos[24 Pi t] + (0.012)Sin[28 Pi t] 
    + (0.011)Cos[28 Pi t] + (0.008)Sin[32 Pi t] + (0.007)Cos[32 Pi t]] \end{verbatim}

\vspace{1cm}

Finally, we can compute the limit $\lim\limits_{E \to -1} \lambda_3^2 (E)$ replacing again the potential in the Sturm-Liouville equation by the first terms of its Fourier series. For $E= - 1$ this amounts to studying  the differential equation

\[ \sin^2 (2 \pi t) A''(t)= \left\{ \frac{1}{2} \pi^2 \Big( 9- 3 \cos (4 \pi t) - 4 \sin (4 \pi t) \Big) - 2 \lambda_3^2 \sin^4 (2 \pi t) \right\} A(t) \]

and the finite-dimensional approximation yields the result

\[ \lim\limits_{E \to - 1} \lambda_3^2 (E) \approx 47.2437 . \]


{\bf Remark:} The second variation formulas prove that, in case of the family \linebreak $g_E =(1+E \cos (2 \pi t))g_o$ $(N=1)$,  the minimal positive eigenvalues of the Laplace and Dirac operator decrease (see Example 4.1) and are smaller than $4 \pi^2$. The numerical evaluation of $\mu_3(E)$ and $\lambda_3^2 (E)$ yields the following table:\\

\begin{center}
\begin{tabular}{|c|c|c|c|c|c|c|c|c|c|} \hline
&&&&&&&&&\\ 
$E$ & 0 & -0.1 & -0.3 & -0.5 & -0.7 & -0.9 & -0.95 & -0.99 & -1\\
&&&&&&&&& \\ \hline  &&&&&&&&&\\
$\mu_3$ & $4 \pi^2$ & 39.284 & 37.897 & 35.741 & 33.378 & 31.09 & 30.5 & 30.1 & 30.013\\
&&&&&&&&&\\ \hline &&&&&&&&&\\
$\lambda_3^2$ & $4 \pi^2$ & 39.333 & 38.353 & 36.714 & 34.983 & 33.331 & 33.2830 & 36.04 & $\approx 36.2$\\
&&&&&&&&&\\ \hline
\end{tabular}
\end{center} 

\vspace{1cm}

\section{Final remarks}

As shown previously, any local deformation $g_E$ of the flat metric realizes the inequality

\[ \mu_1 (g_E) < \lambda_1^2 (g_E) \]

between the first eigenvalues of the Laplace and Dirac operator up to second order. We are not able to give an example of a Riemannian metric $g$ on $T^2$ such that $\lambda_1^2 (g) < \mu_1 (g)$ holds. Moreover, denote again by $\lambda_1^2 (g;l)$ the first positive eigenvalue of the Dirac operator such that the eigenspace contains an $S^1$-representation of weight $l \in {\Bbb Z}$. The corresponding eigenvalue of the Laplace operator we shall denote by $\mu_1 (g;l)$. It is a matter of fact that in all families of Riemannian metrics we have discussed these two eigenvalues are very close. Let us consider, for example, the metric $g_E$ by the function

\[ h_E (t) = e^{\frac{E}{\pi} (\sin (2 \pi t) - 2 \cos (2 \pi t))}    . \]

For the parameter $E=1$ we obtain  the following numerical values using the approximation method described before in the space spanned by the functions $1, \sin (2 \pi nt), \cos (2 \pi nt)$ $(1 \le n \le 5)$:

\[ \lambda_1^2 (g;1)  \approx 6.11056 \quad , \quad \mu_1 (g ;1)  \approx 5.19025  . \]

However, even in this case we already have the inequality $\mu_1 (g_E; 1) < \lambda_1^2 (g_E; 1)$ and the following  figure shows the graph of the two spectral functions for $0 \le E \le 1$ (for the first and the second positive eigenvalue):

\begin{center}

\[
\begin{rotate}[r]{\epsfig{figure=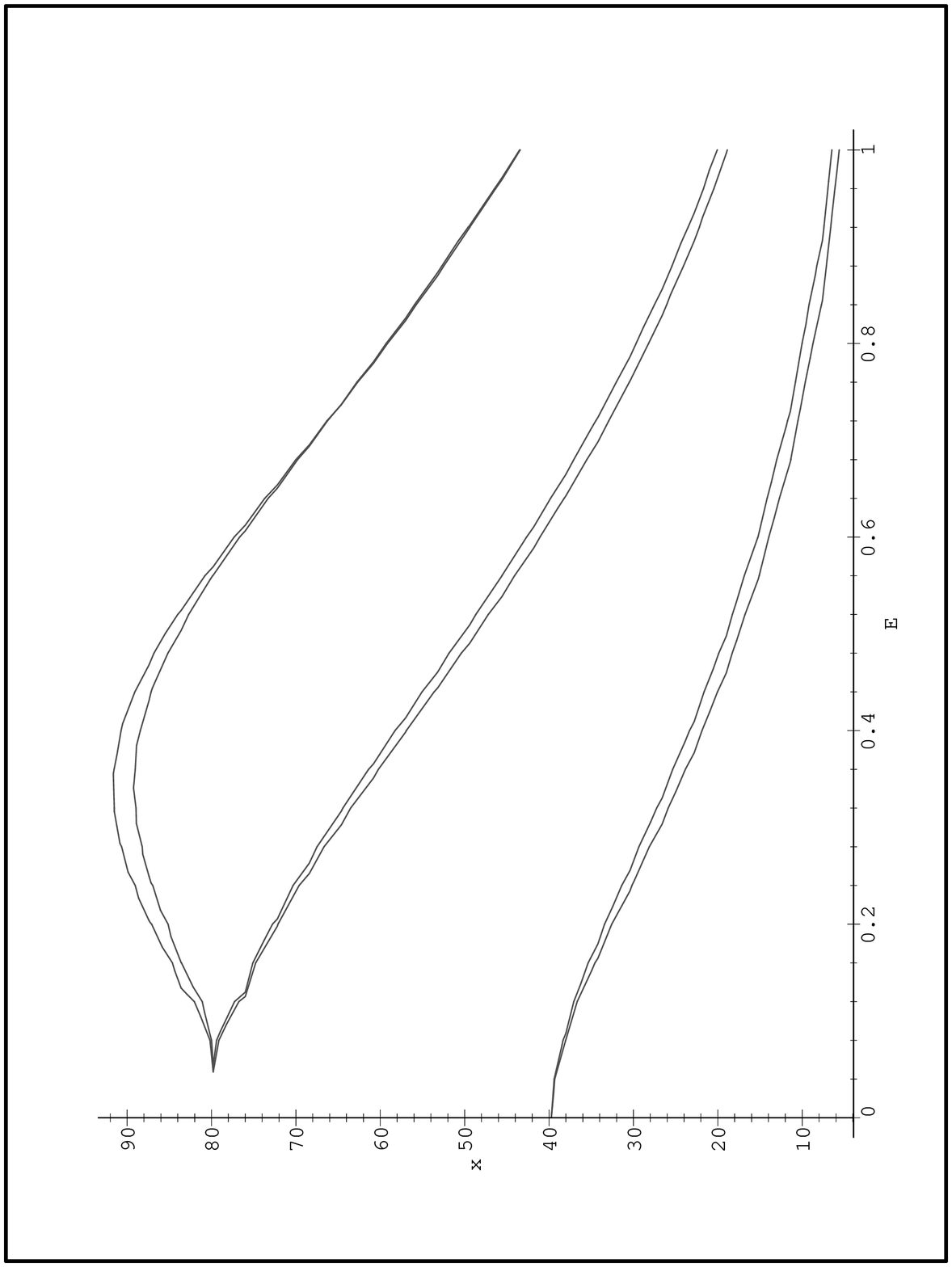,width=10cm}}
\end{rotate}\]

(Figure 10)
\end{center}

\vspace{1cm}


\bigskip


\small
ILKA AGRICOLA\\
Humboldt-Universit\"at zu Berlin, Institut f\"ur Mathematik, Sitz: Ziegelstra\ss e 13a, \\
Unter den Linden 6, D-10099 Berlin\\
{\tt e-mail: agricola@mathematik.hu-berlin.de}\\

BERND AMMANN\\
Universit\"at Freiburg, Mathematisches Institut, Eckerstr. 1, D-79104 Freiburg\\
{\tt e-mail: ammann@mathematik.uni-freiburg.de}\\

THOMAS FRIEDRICH\\
Humboldt-Universit\"at zu Berlin, Institut f\"ur Mathematik, Sitz: Ziegelstra\ss e 13a, \\
Unter den Linden 6, D-10099 Berlin\\
{\tt e-mail: friedric@mathematik.hu-berlin.de}\\

\end{document}